\documentclass[12pt]{amsart}
\usepackage{longtable}
\usepackage{tikz-cd}
\usepackage[total={6.25truein,9truein},centering]{geometry}
\usepackage{hyperref}

\newtheorem{thm}{Theorem}[section]

\newcommand{\CC}{\mathbb C}

\newcommand{\GG}{\mathbb G}

\newcommand{\NN}{\mathbb N}
\newcommand{\ZZ}{\mathbb Z}

\newcommand{\QQ}{\mathbb Q}

\newcommand{\FF}{\mathbb F}
\newcommand{\Cinf}{{\mathbb C}_\infty}

\newcommand{\bh}{\mathbf{h}}

\newcommand{\bm}{\mathbf{m}}

\newcommand{\fp}{\mathfrak{p}}

\newcommand{\tpi}{\widetilde{\pi}}

\DeclareMathOperator{\GL}{GL}
\DeclareMathOperator{\End}{End}
\DeclareMathOperator{\rank}{rank}
\DeclareMathOperator{\Exp}{Exp}
\DeclareMathOperator{\Ext}{Ext}
\DeclareMathOperator{\Gal}{Gal}
\DeclareMathOperator{\Lie}{Lie}
\DeclareMathOperator{\Mat}{Mat}

\DeclareMathOperator{\diag}{diag}

\DeclareMathOperator{\Hom}{Hom}

\newcommand{\tr}{\mathrm{tr}}

\def\as{\mathrel{\raise.095ex\hbox{:}\mkern-4.2mu=}}

\newcommand{\power}[2]{{#1 [\![ #2 ]\!]}}
\newcommand{\laurent}[2]{{#1 (\!( #2 )\!)}}

\title[A rapid introduction to Drinfeld modules, $t$-modules, and $t$-motives]{A rapid introduction to Drinfeld modules, \\ $t$-modules, and $t$-motives}
\author{W.\ Dale Brownawell}
\address{Department of Mathematics, Penn State University,
University Park, PA 16802, USA}
\email{wdb@math.psu.edu}

\author{Matthew A.\ Papanikolas}
\address{Department of Mathematics, Texas A{\&}M University, College
Station, TX 77843, USA}
\email{map@math.tamu.edu}

\date{April 10, 2011; revised May 6, 2016; final version June 11, 2018}

\begin{document}
\maketitle

The theory of Drinfeld modules was initially developed to transport
the classical ideas of lattices and exponential functions to the
function field setting in positive characteristic.  This impulse manifested
itself in work of L.\ Carlitz \cite{Carlitz35}, who
investigated explicit class field theory for $\FF_q(\theta)$ and
defined the later named Carlitz module, which serves as the analogue
of the multiplicative group $\GG_m$.  The theory was brought to rapid
fruition by V.~G.\ Drinfeld \cite{Drinfeld74}, who independently
superseded Carlitz's work and further extended the theory to higher
rank lattices in his investigation of elliptic modules, now called
Drinfeld modules.  G.~W.\ Anderson \cite{Anderson86} saw correctly how
to develop the theory of higher dimensional Drinfeld modules, called
$t$-modules, and at the same time produced a robust motivic
interpretation in his theory of $t$-motives.

The present article aims to provide a brief account of the theories of
Drinfeld modules and Anderson's $t$-modules and $t$-motives.  As such
the article is not meant to be comprehensive, but we have endeavored
to summarize aspects of the theory that are of current interest and to
include a number of examples.  For further information and more complete
details, readers are encouraged to consult the excellent surveys
\cite{DH87,Goss,Hayes92,Rosen,Thakur}.

\section{Exponential functions of algebraic groups}
We begin with some preliminary remarks about commutative algebraic
groups over $\CC$, starting with the multiplicative group $\GG_m$.
The exact sequence
\[
  0 \to 2\pi i \ZZ \to \CC \xrightarrow{\exp} \CC^{\times} \to 1
\]
exhibiting the uniformizability of $\GG_m$,
is the starting point for a multitude of problems and their
solutions in number theory.  For example, the entire study of abelian
extensions of $\QQ$ is intertwined with $\exp(z)$, whose division
values are simply roots of unity.  Moreover, this sequence is the
starting point for transcendence questions involving exponentials
and logarithms of algebraic numbers.

The next natural step along these lines leads to the
investigation of elliptic curves.  We can associate to an elliptic
curve $E$ over $\CC$ a rank $2$ lattice $\Lambda \subseteq \CC$, from which
we can define the Weierstra{\ss} $\wp$-function
\[
  \wp(z) = \frac{1}{z^2} + \sum_{\substack{\lambda \in \Lambda \\ \lambda \neq 0}}
\left( \frac{1}{(z-\lambda)^2} - \frac{1}{\lambda^2} \right).
\]
This leads to the exact sequence
\[
  0 \to \Lambda \to \CC \xrightarrow{\exp_{E}} E(\CC) \to 0,
\]
where $\exp_E(z) = [\wp(z), \wp'(z),1]$, and $\Lambda$ is called the
period lattice of $E$.  When $E$ is defined over a number field $K$,
the division values of $\exp_E(z)$,
much as in the case of $\GG_m$,
generate interesting extensions of $K$ and are the focus of much
study.  Also, transcendence questions about periods and elliptic
logarithms of algebraic points naturally arise from analogy with the
$\GG_m$ case.

These investigations generalize in satisfying ways to general
commutative algebraic groups over $\CC$, including algebraic tori and
abelian varieties.  For a commutative algebraic group $G$ one has the
exponential sequence
\[
  0 \to \Lambda \to \Lie(G) \xrightarrow{\exp_G} G(\CC) \to 0,
\]
where $\Lambda$ is a lattice in $\Lie(G)$, and when $G$ is defined over
a number field, we can similarly study special values of $\exp_G$ and
logarithms of algebraic points on $G$.

Over the past few decades function field analogues have fostered many
fruitful research programs:
\begin{itemize}
\item Cyclotomic theory and explicit class field theory over function fields,
\item Drinfeld modular forms and modular varieties,
\item Drinfeld modules over finite fields,
\item Torsion modules and Galois representations,
\item Characteristic $p$ valued $L$-series,
\item Heights and Drinfeld modules over global function fields,
\item Effective bounds on isogenies of $t$-modules,
\item Transcendence theory,
\item Special functions,
\item Shtukas and automorphic representations over function fields,
\item $t$-motives, $\tau$-sheaves, and Hodge structures.
\end{itemize}
Any complete list of references on the above topics would necessarily
be too long for the scope of this survey.  However, we list here
several useful sources that represent a broad picture and
contain themselves references for further study: texts and monographs
by B\"ockle and Pink~\cite{BocklePink}, Gekeler~\cite{Gekeler},
Goss~\cite{Goss}, Laumon~\cite{Laumon1,Laumon2}, Rosen~\cite{Rosen},
and Thakur~\cite{Thakur}; survey articles by Deligne and
Husemoller \cite{DH87}, Goss~\cite{Goss80a,Goss83,Goss92}, Hartl~\cite{Hartl09},
Hayes~\cite{Hayes92}, Pellarin~\cite{Pellarin08}, and
Thakur~\cite{Thakur95}; and research articles by
Anderson~\cite{Anderson86}, Anderson and Thakur~\cite{AndersonThakur}, David and Denis~\cite{DD99},
Denis~\cite{Denis92}, Drinfeld~\cite{Drinfeld77}, Galovich and
Rosen~\cite{GalRos82}, Gekeler~\cite{Gekeler88}, Gekeler and
Reversat~\cite{GekRev96}, Goss~\cite{Goss79,Goss80b,Goss00},
Hayes~\cite{Hayes85}, Lafforgue~\cite{Lafforgue02}, Pink~\cite{PinkHS,
  Pink97}, Pink and R\"{u}tsche~\cite{PR09}, Poonen~\cite{Poonen95}, Taguchi and Wan~\cite{TagWan96},
Thakur~\cite{Thakur91}, and Yu~\cite{Yu86,Yu91,Yu97}.  Accounts of
many of these topics are also included in the current volume.

\subsection*{Acknowledgements}
This survey was adapted from lecture notes originally written for the
Arizona Winter School at the University of Arizona in 2008.  We thank
the AWS for permitting us to use them here.  We further thank U.~Hartl for making several suggestions that improved the exposition and for pointing out an error in an earlier version.  Research of the second author was supported by NSF Grant DMS-0903838.

\newpage
\section{Drinfeld Modules}
\subsection{Table of  symbols}\

\begin{longtable}{p{0.5truein}@{\hspace{5pt}$\as$\hspace{5pt}}p{4.5truein}}
$p$ & a fixed prime\\
$\FF_q$ & finite field of $q = p^m$ elements\\
$\FF_q[t]$ & polynomials in the variable $t$\\
$k$ & $\FF_q(\theta)=$ rational functions in the variable $\theta$\\
$k_{\infty}$ & $\FF_q(\!(1/\theta)\!)=$ $\infty$-adic completion of $k$\\
$|\cdot|_\infty$ & absolute value on $k_\infty$ such that $|\theta|_\infty=q$\\
$\overline{k_{\infty}}$ & algebraic closure of
$k_{\infty}$\\
$\CC_\infty$ & completion of $\overline{k_{\infty}}$ with respect
to $|\cdot|_\infty$ \\
$\overline{k}$ & algebraic closure of $k$ in $\CC_\infty$ \\
$\deg$ & the function associating to each element of $\FF_q[\theta]$
its degree in $\theta$\\
$\tau$ & the $q$-power Frobenius map sending $x \mapsto x^q$ on a commutative $\FF_q$-algebra $R$\\
$c^{(i)}$ & $c^{q^i}=$ the $i$th iterate of $\tau$ applied to an element $c \in R$, $i \in \ZZ$\\
$R\{\tau\}$ & the ring of \emph{twisted polynomials} $\sum_i a_i\tau^i$, $a_i \in R$, where multiplication is given by $a\tau^i \cdot b \tau^j = ab^{(i)}\tau^{i+j}$\\
$R\{\sigma\}$ & the ring of twisted polynomials $\sum_i a_i\sigma^i$, $a_i \in R$, $R$ perfect, where multiplication is given by $a\sigma^i \cdot b \sigma^j = ab^{(-i)}\sigma^{i+j}$
\end{longtable}

\subsection{The Carlitz module}
The Carlitz module $C$ is the first example of a Drinfeld module.
Defined by Carlitz \cite{Carlitz35} in 1935, it is given by the
$\FF_q$-algebra homomorphism
\[
  C : \FF_q[t] \to \CC_\infty \{ \tau \}
\]
defined so that
\[
C(t) = \theta + \tau.
\]
(Truth be told, Carlitz set $C(t) = \theta-\tau$, but the definition
we are using is more prevalent today.)  The natural point of view is that a twisted polynomial $f = a_0 + a_1\tau + \cdots + a_d \tau^d \in
\CC_\infty \{ \tau \}$ represents the $\FF_q$-linear endomorphism of
$\CC_\infty$,
\[
  x \mapsto f(x) = a_0x + a_1x^q + \cdots + a_d x^{q^d}.
\]
In this way $C$ makes $\CC_\infty$ into an $\FF_q[t]$-module, where $a \cdot
x = C(a)(x)$, and in particular
\[
  C(t)(x) = \theta x + x^q.
\]

Exponential functions enter the picture with the Carlitz exponential
function
\[
  \exp_C(z) \as \sum_{i \geq 0} \frac{z^{q^i}}{D_i},
\]
where $D_0 = 1$ and $D_i = (\theta^{q^i} - \theta)(\theta^{q^i} -
\theta^q) \cdots (\theta^{q^i} - \theta^{q^{i-1}})$ for $i \geq 1$.
This function converges for all $z \in \CC_\infty$, and the recursion
$D_i = (\theta^{q^i}-\theta)D_{i-1}^q$ implies that
\[
  \exp_C(\theta z) = \theta\exp_C(z) + \exp_C(z)^q = C(t)(\exp_C(z)).
\]
More generally, one checks that the following diagram commutes for any $a \in \FF_q[t]$:
\[
\begin{tikzcd}
\CC_\infty \ar[r,"\exp_C(z)"] \ar[d,"z\,\mapsto\,a(\theta)z"']
& \CC_\infty \ar[d,"x\,\mapsto\,C(a)(x)"] \\
\CC_\infty \ar[r,"\exp_C(z)"'] & \CC_\infty.
\end{tikzcd}
\]
Thus $\exp_C : \Cinf \to \Cinf$ is an $\FF_q[t]$-module
homomorphism, where $t$ acts on the domain by scalar multiplication by
$\theta$ and on the range by the endomorphism $C(t)$.  When convenient
we will use $(C,\Cinf)$ to denote $\CC_\infty$ with the Carlitz
$\FF_q[t]$-module structure.

The Carlitz exponential function uniformizes the Carlitz module
as follows.  As $\CC_\infty$ is algebraically closed, it
follows from the Weierstrass preparation theorem that the Carlitz
exponential is surjective.  Remarkably, Carlitz found its kernel to
be all $\FF_q[\theta]$-multiples of
\begin{equation} \label{CarlitzPer}
\tpi \as \theta \sqrt[q-1]{-\theta} \prod_{i=1}^\infty \Bigl( 1- \theta^{1-q^i}
\Bigr)^{-1} \in k_\infty \bigl( \sqrt[q-1]{-\theta} \bigr),
\end{equation}
where $\sqrt[q-1]{-\theta}$ is any fixed root of $-\theta$ (see \cite[\S 3.2]{Goss}, \cite[\S 2.5]{Thakur}.  We then
have an exact sequence of $\FF_q[t]$-modules,
\[
  0 \to \FF_q[\theta] \cdot \tpi \to \CC_\infty
\xrightarrow{\exp_C} (C,\Cinf) \to 0.
\]
The quantity $\tpi$ is called the \emph{Carlitz period}.  This uniformization sequence underlies the first transcendence results in positive characteristic: in 1941, L.~I.\ Wade \cite{Wade} showed that $\tpi$ is transcendental over $\FF_q(\theta)$.

One of Carlitz's motivations for studying the Carlitz module was to
explore explicit class field theory for the rational function field
$\FF_q(\theta)$, as in \cite{Carlitz38}.  For $f \in
\FF_q[t]$, we let
\[
  C[f] \as \{ x \in \CC_\infty \mid C(f)(x) = 0 \}
\]
denote the \emph{$f$-torsion} on $C$, which is isomorphic to
$\FF_q[t]/(f)$ as an $\FF_q[t]$-module and which is a Galois module over the
separable closure of $\FF_q(\theta)$.  The Carlitz cyclotomic field is
the field $\FF_q(\theta,C[f])$, and there is an isomorphism
\[
  \rho : (\FF_q[t]/(f))^{\times} \stackrel{\sim}{\longrightarrow} \Gal( \FF_q(\theta,C[f])/ \FF_q(\theta) ),
\]
such that for an $\FF_q[t]$-module generator $\lambda \in C[f]$ we
have $\rho_a(\lambda) = C(a)(\lambda)$ for any $a \in
(\FF_q[t]/(f))^{\times}$.  Moreover, $\rho_a$ coincides with the Artin
automorphism for $a$, and in this way we obtain an explicit Galois
action on a piece of the maximal abelian extension of $\FF_q(\theta)$
that agrees with class field theory.  However, at $\infty$ Carlitz's cyclotomic extensions
are at most tamely ramified.  Only later did Hayes
\cite{Hayes74} complete the picture to obtain a full analogue of the
Kronecker-Weber Theorem, by showing how Carlitz's constructions could
be used to describe the Galois action on the full maximal abelian
extension of $\FF_q(\theta)$ and provide an explicit class field
theory as well.  Indeed, the similarities we observe between the
Carlitz module and the multiplicative group identifies a theme that
pervades the theory: in the dictionary between function fields and
number fields we have
\[
  C \longleftrightarrow \GG_m,
\]
and this identification often occurs even when not completely anticipated.

\subsection{Drinfeld modules}
After Carlitz the situation became clearer through the work of V.~G.\
Drinfeld~\cite{Drinfeld74,Drinfeld77} and D.~Hayes
\cite{Hayes74,Hayes79} in the 1970's.  Drinfeld introduced what he
called \emph{elliptic modules} (now commonly called \emph{Drinfeld
  modules}) because they have remarkable similarities with classical
elliptic curves.  Drinfeld simultaneously generalized Carlitz's work
in two directions: he extended the definitions to arbitrary rings of
functions on curves over finite fields and to arbitrary rank lattices.
While Drinfeld was unaware of Carlitz's previous work, Hayes
continued Carlitz's work on explicit class field theory of arbitrary
function fields of one variable over finite fields and developed a
rank $1$ theory that coincided with Drinfeld's.  (Consequently rank
$1$ ``sign normalized'' Drinfeld modules are often referred to as ``Drinfeld-Hayes''
modules.)  It is worth pointing out that Carlitz had conceived of a
theory of exponential functions for lattices of rank higher than $1$,
as evidenced by his manuscript \cite{Carlitz95} not published until
1995.

The essential construction of a Drinfeld module is the following.  Let
$L$ be an arbitrary extension of $\FF_q$.  Let $A$ be the ring of
functions on a smooth projective geometrically irreducible curve $X/\FF_q$ that are regular away
from a fixed $\FF_q$-rational point $\infty \in X$, and fix an
$\FF_q$-algebra homomorphism $\iota: A \to L$.  In general $\infty$ need not be $\FF_q$-rational, but we consider this case here for convenience.  Readers new to the
subject are encouraged to take simply $A = \FF_q[t]$, for which the
theory is just as rich as in the general case but at times more
straightforward.  A Drinfeld $A$-module is then an $\FF_q$-algebra
homomorphism
\[
\phi \colon A \rightarrow L \{\tau\}
\]
for which $\phi(a) = \iota(a)\tau^0 +{}$ higher order terms in $\tau$.
As in the case of the Carlitz module, elements of the ring
$L\{\tau\}$ can be thought of as $\FF_q$-linear endomorphisms of the additive group of
$L$, and thus we often identify $\phi$ with the $A$-module structure
on $L$ induced by $\phi$ and write $(\phi,L)$ for $L$ with this new
$A$-module structure.  Geometrically, the Drinfeld module is simply
the additive group $\GG_a$ over $L$, but we think of $\phi$ as the
pair $(\phi,\GG_a)$, where
\[
  \phi \colon A \to \End_L(\GG_a),
\]
and $\End_L(\GG_a)$ consists of all endomorphisms defined over $L$.

If $\iota : A \to L$ is injective, then $\phi$ is said to have
\emph{generic characteristic}.  If not, then $\phi$ is said to have
\emph{characteristic $\fp$}, where $\fp = \ker \iota \neq (0)$. If
$\phi(A) \subseteq K \{\tau\}$ for some subfield $K \subseteq L$, we
say that $\phi$ is \emph{defined over $K$}.  There is a non-negative
integer $r$ such that for every $a \in A$,
\[
  \deg_{\tau}(\phi(a)) = r \deg(a),
\]
where $\deg(a)$ is normalized by $|A/(a)| = q^{\deg(a)}$ (see \cite[\S 4.5]{Goss}).  The integer
$r$ is called the \emph{rank} of $\phi$.  Although it is not
immediate, we will see in the next section that there exist Drinfeld
modules of any rank $r$ for any ring $A$.  In the case of the Carlitz
module we have $A = \FF_q[t]$, the structure morphism $\iota : A \to
\CC_\infty$ is defined by $\iota(a) = a(\theta)$, and the Carlitz
module is a rank $1$ Drinfeld $A$-module with generic characteristic.

\subsection{Drinfeld modules and lattices} \label{DrinfeldMods}
We continue with the notation of the previous section.  In the case
that we have an embedding $\iota : A \hookrightarrow \CC_\infty$,
Drinfeld constructed Drinfeld $A$-modules in close analogy to the
situation of elliptic curves and elliptic functions over $\CC$.  We
consider $\CC_\infty$ to be an $A$-module via $\iota$, and
for this section we identify $A$
with its image $\iota(A)$ in $\CC_\infty$.

Starting with an \emph{$A$-lattice} $\Lambda \subseteq \CC_\infty$,
i.e.\ a discrete finitely generated projective $A$-submodule of
$\CC_\infty$, of rank $r > 0$, we define the lattice function
\[
  \exp_\Lambda(z) := z {\prod_{\lambda \in \Lambda}}' \biggl( 1 - \frac{z}{\lambda} \biggr), \quad z \in \CC_\infty,
\]
where the product is taken over all non-zero lattice elements.  The
discreteness of $\Lambda$ ensures that $\exp_\Lambda(z)$ converges for
all $z \in \CC_\infty$: only finitely many $\lambda \in \Lambda$ lie
within any given distance to the origin.  Consideration of the partial
products involving bounded $\lambda$ shows that the Drinfeld
exponential function has an expansion of the form
\[
\exp_\Lambda(z) = z + \sum_{i \geq 1} a_i z^{q^i},
\]
and as such it is an \emph{$\FF_q$-linear} power series.  Thus,
for $c \in \FF_q$,
\[
\exp_\Lambda(z_1 + c z_2) = \exp(z_1) + c\exp_\Lambda(z_2).
\]
Also for $c \in \CC_\infty$ non-zero, the product expansion for
$\exp_\Lambda(z)$ makes obvious that
\[
\exp_{c\Lambda}(cz) = c\exp_\Lambda(z).
\]
Moreover $\exp_\Lambda(z)$ visibly parametrizes $\CC_\infty$ and has
kernel $\Lambda$.  So the sequence
\[
0 \to \Lambda \to \CC_\infty \xrightarrow{\exp_\Lambda} \CC_\infty \to 0
\]
is exact.

Now consider the case that $\Lambda_1$, $\Lambda_2$ are $A$-lattices
of the same $A$-rank, but $\Lambda_1 \subseteq \Lambda_2$.  Then
$\Lambda_2/\Lambda_1$ is a finite dimensional $\FF_q$-vector space,
say with coset representatives $\lambda_0(= 0), \dots, \lambda_{d-1}$.
Hence
\[
P_{[\Lambda_2\colon \Lambda_1]}(X) \as X \prod_{i=1}^{d-1}
 \left( 1 -  \frac{X}{\exp_{\Lambda_1}(\lambda_i)}\right)
\]
is an $\FF_q$-linear polynomial in $X$ with $X$ as lowest term, and it
provides the crucial functional relation
\[
\exp_{\Lambda_2}(z) = P_{[\Lambda_2\colon \Lambda_1]}(\exp_{\Lambda_1}(z)),
\]
since both sides have the same simple zeros and the same leading terms.

When $\Lambda_1$, $\Lambda_2$ have the same rank and
$c\Lambda_1 \subseteq \Lambda_2$, then $[\Lambda_2 \colon c\Lambda_1]$
is finite.  In that case,
\[
\exp_{\Lambda_2}(cz) = P_{[\Lambda_2 \colon c\Lambda_1]}(\exp_{c\Lambda_1}(cz))
= P_{[\Lambda_2 \colon c\Lambda_1]}(c\exp_{\Lambda_1}(z))
\]
as both sides have the same zeros and the same leading terms.
In particular, when $\Lambda_1 = \Lambda_2 = \Lambda$ and $c =
a \in A$, $a \neq 0$, we can write
\begin{equation} \label{FunctionalEqn}
\exp_\Lambda(az) = \phi_\Lambda(a) \exp_\Lambda(z)
\end{equation}
for $\phi_\Lambda(a)(z) = P_{[\Lambda \colon a\Lambda]}(z)$, where $
\phi_\Lambda(a) = a\tau^0 +{}$ higher order terms in $\tau$ lies in
$\CC_\infty\{\tau\}$ and is non-zero.  Thus we have defined a
function
\begin{equation}
  \phi_\Lambda : A \to \CC_\infty\{ \tau \},
\end{equation}
and we will see shortly that $\phi_\Lambda$ defines a Drinfeld $A$-module.

Returning to the case $c\Lambda_1 \subseteq
\Lambda_2$, $c \in \CC_\infty$, and writing $P_{[\Lambda_2 \colon
  c\Lambda_1]}(cz) = \psi_{(c)}(z)$ for $\psi_{(c)}$ in
$\CC_\infty\{\tau\}$, we find
\begin{multline*}
  \phi_{\Lambda_2}(a)\psi_{(c)} \exp_{\Lambda_1}(z) = \phi_{\Lambda_2}(a)
  \exp_{\Lambda_2}(cz) = \exp_{\Lambda_2}(acz) \\
  = \psi_{(c)} \exp_{\Lambda_1}(az) = \psi_{(c)}
  \phi_{\Lambda_1}(a)\exp_{\Lambda_1}(z).
\end{multline*}
Since $\exp_{\Lambda_2}(z)$ is a transcendental function (it has
infinitely many zeros), we conclude that
\[
\phi_{\Lambda_2}(a)\psi_{(c)} = \psi_{(c)} \phi_{\Lambda_1}(a),
\]
for all $a \in A$, and we say that $\psi_{(c)} \in \Hom_A(\phi_{\Lambda_1},\phi_{\Lambda_2})$.  Any non-zero element $\psi$ of
$\CC_\infty\{\tau\}$ satisfying this property will be called an
\emph{isogeny} from $\phi_{\Lambda_1}$ to $\phi_{\Lambda_2}$, and we also
write
\[
\psi \colon \phi_{\Lambda_1} \rightarrow \phi_{\Lambda_2}.
\]
Isogeneity is an equivalence relation.  Note that, as $c \ne 0$,
$\psi_{(c)} = c\tau^0 +{}$ higher order terms in $\tau$ in the previous
displayed line.  It is not hard to see that if such an isogeny has the form $\psi=
c\tau^0 + {}$ higher order terms, then $\psi = \psi_{(c)}$.

In particular, when $ \Lambda_1 = \Lambda_2 = \Lambda$ and $c \in A$,
we have that $\psi_{(c)} = \phi_{\Lambda}(c)$ and
\[
\phi_{\Lambda}(a) \phi_\Lambda(c) \exp_\Lambda(z) = \exp_\Lambda (acz)
= \phi_{\Lambda}(ac) \exp_\Lambda(z).
\]
Since $a \mapsto \phi_\Lambda (a)$ is also additive, $\phi_\Lambda$ is
thus a ring homomorphism, and we have a Drinfeld module structure
induced by the effect on the range of $\exp_\Lambda$, i.e.\ on $\Cinf$,
of the $A$-action on $\Lambda$.  Furthermore, since the degree in $z$
of $\phi_\Lambda(a)(z)$ is $[\Lambda:a\Lambda] = q^{r\deg{a}}$, we see
that the degree in $\tau$ of $\phi_\Lambda$ is $r$.  Thus
$\phi_\Lambda$ is a Drinfeld $A$-module of rank $r$.

One also checks, for $\Lambda_1$, $\Lambda_2$ as above, that
$\Hom_A(\phi_{\Lambda_1},\phi_{\Lambda_2})$ is an
$A$-module, and so for any particular lattice $\Lambda$ the
endomorphism ring of $\phi_{\Lambda}$ is the $A$-algebra
$\End_A(\phi_\Lambda) \as \Hom_A(\phi_\Lambda,\phi_\Lambda)$.  Since
$\End_A(\phi_\Lambda)$ can be identified with those $c \in \Cinf$ such
that $c\Lambda \subseteq \Lambda$, it follows that
$\End_A(\phi_\Lambda)$ is an integral domain and is a finitely
generated $A$-module of projective rank at most $r$ \cite[Ch.~4]{Goss}.

In generic characteristic Drinfeld demonstrated the striking fact that the considerations proceeding from the analytic to the algebraic are also reversible.

\begin{thm}[Drinfeld's Uniformization Theorem] Given a homomorphism
  $\phi : A \rightarrow \CC_\infty\{\tau\}$ such that $\phi(a) =
  a + \dots + a_m \tau^m$, $a_m \ne 0$, $m=r\deg(a)$, there is a
  unique $A$-lattice $\Lambda$ such that $\phi = \phi_\Lambda$.
  Moreover $\rank_A \Lambda = r$.
\end{thm}

A brief outline of how to see this is to use the condition $e(\theta
z)=\phi(t)e(z)$ to define a unique $\FF_q$-linear power series
$e_{\phi}(z)$ with leading term $z$.  Then one shows that
$e_{\phi}(z)$ is $\FF_q$-linear and entire.  Finally, from the
functional equation, one sees that the zeros of $e_{\phi}(z)$ form a
discrete $A$-module, i.e.\ a lattice $\Lambda$, and then it remains to
show that $e_{\phi}(z) = \exp_\Lambda(z)$ by the uniqueness of the solution to \eqref{FunctionalEqn} having lowest term $z$.  For a complete proof, see
\cite{Goss, Rosen, Thakur}.

\subsection{The Weierstra{\ss}-Drinfeld correspondence}
One is reminded of the situation of elliptic curves over the
complex numbers, and the analogies are amazingly tight.  Based on our
various observations, we have the following dictionary.

\bigskip
\begin{center}
\begin{tabular}{p{2.75in}@{\hspace{0.25truein}}p{2.75in}}
  Weierstra{\ss} & Drinfeld \\ \hline \\
  $\ZZ$ & $A$ \\
  $2$-dim.\ lattice $\Lambda$ & $r$-dim.\ $A$-lattice $\Lambda$ \\
  $\exp_\Lambda(z) = (\wp(z),\wp^\prime(z))$ &
  $\exp_{\Lambda}(z) \as z\, {\prod}'_{\lambda \in \Lambda}( 1 - \frac{z}{\lambda})$ analytic!\\
  elliptic curve $\mathcal{E} \colon
    y^2 = 4x^3 -g_2 x - g_3$ &  $\GG_a$ \\
$0 \rightarrow \Lambda \rightarrow \CC \rightarrow \mathcal{E}_\Lambda(\CC)
\rightarrow 0$ &
$0 \rightarrow \Lambda \rightarrow \CC_\infty \rightarrow (\phi_\Lambda,\Cinf)
\rightarrow 0$ \\
Isogenies given as $c$ s.t.\ $c\Lambda_1 \subseteq \Lambda_2$ &
Isogenies given as $c$ s.t.\ $c\Lambda_1 \subseteq \Lambda_2$\\
& \hspace*{10pt}when $\rank_A\Lambda_1 = \rank_A \Lambda_2$ \\
$\ZZ \subseteq \End(\mathcal{E})$  & $\phi_\Lambda \colon A \rightarrow \End(\GG_a)$ via $\phi_\Lambda(t) = \theta \tau^0 +{}$ \\
& \hspace*{10pt}higher terms $\in \CC_\infty\{\tau\}$
\end{tabular}
\end{center}

\section{\texorpdfstring{$t$}{t}-Modules}
By analogy with taking the step from elliptic curves to abelian
varieties, one can ask questions about how to define higher
dimensional Drinfeld modules properly.  In 1986, Anderson \cite{Anderson86}
devised and solved this problem by defining \emph{$t$-modules}.
Moreover, his construction intrinsically includes many reasonable
generalizations of Drinfeld modules to the higher dimensional setting,
including direct products, tensor products, and extensions.  Anderson
also defined a category of companion objects called
\emph{$t$-motives}, which will be the subject of \S\ref{tMotives}.

\subsection{Definitions}
Throughout this section we assume that $A = \FF_q[t]$.  It is possible
to define a theory of `$A$-modules' and `$A$-motives' for general $A$,
but to simplify things we adhere to the $\FF_q[t]$ case as in
\cite{Anderson86}.  For a commutative $\FF_q$-algebra $R$, a matrix
$B \in \Mat_{m\times n}(R)$, and $i \in \NN \cup \{0\}$, we set $B^{(i)}$  to be
the matrix whose $jk$-entry is $B_{jk}^{(i)} = B_{jk}^{q^i}$.  In this way, when
$m=n$, the map $B \mapsto B^{(i)}$ is an $\FF_q$-algebra homomorphism
$\Mat_n(R) \to \Mat_n(R)$ and we can define the ring of twisted
polynomials $\Mat_n(R)\{\tau\}$ so that
\[
  B\tau^i \cdot C \tau^j = BC^{(i)}\tau^{i+j}.
\]
Thus we really have $\Mat_n(R)\{\tau\} = \Mat_n(R\{\tau\})$.  Moreover, we can map $\Mat_n(R)\{\tau\}$ to a subring of $\End_R(\GG_a^n)$, where for
$\sum B_i\tau^i \in \Mat_n(R)\{\tau\}$ and $x \in \Mat_{n\times 1}(R)$
\begin{equation} \label{matrixtwists}
  \left( \sum B_i \tau^i \right) (x) = \sum B_i x^{(i)}.
\end{equation}

Now let $L$ be an extension of $\FF_q$, and fix $\iota \colon A \to L$.
An \emph{Anderson $t$-module} over $L$ is then defined by an $\FF_q$-algebra
homomorphism
\[
  \Phi \colon A \to \Mat_d(L)\{\tau\},
\]
such that if we set
\[
  \Phi(a) := \partial\Phi(a)\tau^0 + {}\textnormal{higher order terms in $\tau$},
\]
where $\partial$ denotes the differentiation map on $\GG_a^d$ at the origin,
then $\partial\Phi(t) = \theta I_d + N$, where~$N$ is a nilpotent matrix and
$I_d$ is the $d\times d$ identity matrix.  In this way $A$ operates on
$\Mat_{d\times 1}(L) = \GG_a^d(L)$ via $\Phi$ through
\eqref{matrixtwists}, and we will often say that a $t$-module is given
by the pair $(\Phi,\GG_a^d)$ or $(\Phi,L^d)$ to denote this action.  We say that
$\Phi$ has \emph{dimension $d$}, and thus a $1$-dimensional $t$-module is
simply a Drinfeld $A$-module.

If $L = \Cinf$, then we can also define a unique exponential function
$\Exp_\Phi \colon \Cinf^d \to \Cinf^d$, via a power series in
$z_1, \dots, z_d$,
\[
\Exp_\Phi(z) = z + \sum_{i \geq 1} B_i z^{(i)}, \quad z =
\left( \begin{smallmatrix} z_1 \\ \vdots \\ z_d \end{smallmatrix}
\right),\ B_i \in \Mat_d(\Cinf),
\]
satisfying for $a = t$ and thus for all $a \in A$,
\[
  \Exp_\Phi(\partial\Phi(a)z) = \Phi(a)\left( \Exp_\Phi(z) \right).
\]
This functional equation uniquely determines the coefficients $B_i$, $i \geq 1$.
The function $\Exp_\Phi$ converges on all of $\Cinf^d$, and if
$\Exp_\Phi$ is surjective, then $E$ is said to be
\emph{uniformizable}.  Surjectivity of the exponential map is somewhat
subtle (see \cite[{\S 2.2}]{Anderson86}) and is not guaranteed.
However, all exponential maps occurring in this note are surjective;
in other words all $t$-modules we will consider are uniformizable.

The kernel of $\Exp_\Phi$ is a $\partial\Phi(A)$-submodule $\Lambda$ of
$\Cinf^d$, which is finitely generated and discrete.  Just as in the
case of abelian varieties, not every $\partial\Phi(A)$-lattice in $\Cinf^d$
is the kernel of an exponential function for some $t$-module.  We define the \emph{rank} of $\Phi$ to be
the rank of $\Lambda$ as a $\partial\Phi(A)$-module.  Thus if $E$ is
uniformizable, the exponential function induces a familiar exact
sequence of $A$-modules
\[
  0 \to \Lambda \to \Cinf^d \xrightarrow{\Exp_\Phi} (\Phi,\Cinf^d) \to 0
\]
For proofs of the above statements about $\Exp_\Phi$, see
\cite{Anderson86,Goss}.

We summarize the connections between Drinfeld modules and $t$-modules:

\bigskip
\begin{center}
\begin{tabular}{p{2.75in}@{\hspace{0.25truein}}p{2.75in}}
  Drinfeld $A$-modules & Anderson $t$-modules
\\ \hline \\
$\phi(t) \in \CC_\infty\{\tau\}$  &  $\Phi(t) \in \Mat_d(\CC_\infty)\{\tau\}$ \\
$\phi(t) = \theta \tau^0 +{}$ higher order terms &
$\Phi(t) = \partial\Phi(t)\tau^0 +{}$ higher order terms \\
$\partial\phi(t) = \theta$  & $\partial\Phi(t) = \theta I_d + N \in \Mat_d(\Cinf)$\\
  --- & $N^d = 0$ \\
unique entire $\exp_\phi \colon \CC_\infty \rightarrow \CC_\infty$ &
unique entire $\Exp_\Phi \colon \CC_\infty^d \rightarrow
\CC_\infty^d$\\
$\exp_\phi(\theta z) = \phi(t)(\exp_\Lambda(z))$ &
$\Exp_\Phi(\partial\Phi(t) \mathbf{z}) = \Phi(t)\Exp_\Phi(z)$\\
$\exp_\phi(z) = z + \sum_{i \geq 1} b_i z^{q^i}$  &
$\Exp_\Phi(z) = z + \sum_{i \geq 1} B_i z^{(i)}$ \\
$\Lambda = \ker(\exp_\phi)$  finitely generated &
$\Lambda = \ker(\Exp_\Phi)$  finitely generated \\
\hspace*{10pt}discrete $A$-submodule of $\CC_\infty$ & \hspace*{10pt}discrete
$\partial\Phi(A)$-submodule of $\CC_\infty^d$ \\
$\exp_\phi$ always surjective on $\Cinf$ & surjectivity of $\Exp_\Phi$ not guaranteed
\end{tabular}
\end{center}

\bigskip Although the functional equation of the exponential function
has a unique solution, as we have noted above, it is perhaps of
passing interest that two different $t$-modules may have the same
exponential function \cite{Brown01}.

Let $E = (\Phi,\GG_a^d)$ and $F = (\Psi,\GG_a^m)$ denote two
$t$-modules over a field $L$.  Then by a \emph{morphism} $f \colon F
\to E$ over $L$, we mean a morphism of commutative algebraic groups $f
\colon \GG_a^m \to \GG_a^d$ over $L$ commuting with the action of $A$:
\[
  f \Psi(t) = \Phi(t)f.
\]
A \emph{sub-$t$-module} of $E$ is then defined to be the image of any closed immersion
$f \colon F \to E$, which is itself isomorphic as an algebraic group
to $\GG_a^s$ for some $s$, is invariant under the $A$-action, and is
isomorphic to a $t$-module.

When $L = \Cinf$, we can further describe sub-$t$-modules of a
$t$-module $E = (\Phi,\GG_a^d)$ as follows.  By identifying both
$\Lie(E)$ (the tangent space at the origin of $\GG_a^d$) and $E$ with
copies of $\Cinf^d$, we have
\[
  \Exp_\Phi \colon \Lie(E) \to E.
\]
In this setting, a sub-$t$-module is a connected algebraic
subgroup $F$ of $\GG_a^d$ such that
\begin{itemize}
\item $\partial\Phi(t)(\Lie(F)) \subseteq \Lie(F)$,
\item $\Exp_\Phi(\Lie(F)) = F(\CC_\infty)$.
\end{itemize}
In other words, $t$-modules satisfy the usual Lie correspondence for
algebraic groups.

Now we turn to some interesting examples of $t$-modules:

\subsection{Products of other \texorpdfstring{$t$}{t}-modules}
In particular, if $\phi_{1}, \dots, \phi_{n}$ are Drinfeld modules,
then taking $\Phi(t) = \diag(\phi_{1}(t),\dots,\phi_{n}(t))$ sets $N =
0$ and gives rise to $\Exp_{\Phi}(z) =
(\exp_{\phi_{1}}(z_{1}),\dots,\exp_{\phi_{n}}(z_{n}))^{\tr}$.

\subsection{Tensor powers of the Carlitz module} \label{CarlitzTensorMod}
The tensor powers $C^{\otimes n}$, $n \geq 1$, of the Carlitz module
were defined and investigated extensively by Anderson and Thakur in
\cite{AndersonThakur}.  That $C^{\otimes n}$ is the $n$-fold tensor
product of $C$ relies on the tensor product construction in
\S\ref{tMotives}, but we can define these $t$-modules directly and
extract many interesting properties without this information.

We define
\[
  C^{\otimes n} \colon A \to \Mat_n(\Cinf)\{\tau\}
\]
by
\[
  C^{\otimes n}(t) = \theta I_n + N + E\tau,
\]
where
\[
  N = \begin{pmatrix} 0 & 1 & \cdots & 0 \\ \vdots & \ddots  & \ddots
& \vdots \\
\vdots & &  \ddots & 1 \\ 0 & \cdots & \cdots & 0 \end{pmatrix}, \qquad
E = \begin{pmatrix} 0 & \cdots & \cdots & 0\\ \vdots & & & \vdots \\
\vdots & & & \vdots \\ 1 & \cdots & \cdots & 0 \end{pmatrix}.
\]
Of particular interest are the following two results obtained in
\cite{AndersonThakur}.  First, if we let
\[
  \Exp_n \colon \Cinf^n \to \Cinf^n
\]
be the exponential function of $C^{\otimes n}$, then there is a vector
\[
  \pi_n = \begin{pmatrix} * \\ \vdots \\ * \\ \tpi^n \end{pmatrix} \in
\Cinf^n,
\]
such that
\[
  \ker(\Exp_n) = dC^{\otimes n}(A) \cdot \pi_n.
\]
Thus $C^{\otimes n}$ has rank $1$, and the $n$-th power of $\tpi$ is
the final coordinate of a generator of the period lattice for
$C^{\otimes n}$.

Anderson and Thakur also show that the Carlitz zeta value
\[
  \zeta_C(n) := \sum_{\substack{a \in A \\ \textnormal{monic}}} \frac{1}{a^n}
\in k_\infty,
\]
is closely involved with $C^{\otimes n}$.  Moreover, they find
explicit points
\[
  s_n = \begin{pmatrix} * \\ \vdots \\ * \\ \Gamma_n \zeta_C(n)
\end{pmatrix} \in \Cinf^n, \quad S_n \in (C^{\otimes n},\FF_q(\theta)^n),
\]
such that
\[
  \Exp_n(s_n) = S_n.
\]
Here $\Gamma_n \in A$ is the Carlitz factorial \cite{Goss,Thakur}.  Thus, $\zeta_C(n)$ is
the coordinate of the logarithm of a point on $C^{\otimes n}$ that is
defined over $\FF_q(\theta)$.  Furthermore they prove that $S_n$ is a
torsion point on $C^{\otimes n}$ if and only if $(q-1)\mid n$, which
is intertwined with the result of Carlitz that
\[
  (q-1) \mid n \Rightarrow \frac{\zeta_C(n)}{\tpi^n} \in
\FF_q(\theta)^\times.
\]
See \cite{AndersonThakur} for more details.

\subsection{\texorpdfstring{$t$}{t}-modules arising from quasi-periodic functions} \label{QuasiPerMods}
If one is led by Drinfeld to natural function field analogies with elliptic
curves, one can also be inspired to pursue a further analogy with the
elliptic situation---that of extensions of elliptic curves
$\mathcal{E}$ by the additive group $\GG_a$, giving rise to
quasi-elliptic functions:
\[
0 \rightarrow \GG_a \rightarrow E \rightarrow \mathcal{E} \rightarrow 0,
\]
where the exponential function of $E$ is given by
\[
(z,u) \longmapsto (1,\wp(z), \wp^\prime(z), u - \zeta(z))
\]
and $\zeta(z)$ is the quasi-periodic Weierstra{\ss} zeta function.  Here $z$ is the coordinate on $\Lie(\mathcal{E})$ and $u$ is the coordinate on $\Lie(\GG_a)$; the first three coordinates in the image of this formula are the projective coordinates on $\mathcal{E}$ and the fourth is the affine coordinate on $\GG_a$.
The periods of this map are the pairs $(\omega,\eta)$, where $\omega =
n_1\omega_1 + n_2\omega_2$ is a period of $\wp(z)$ and $\eta = n_1
\eta_1 + n_2 \eta_2$ is the corresponding quasi-period expressed in
terms of a basis $\omega_1, \omega_2$ of periods for $\wp(z)$ and
$\eta_i = 2\zeta(\omega_i/2)$, $i = 1,2$.

Anderson, Deligne, Gekeler, and Yu developed a theory of
quasi-periodic Drinfeld functions \cite{Gekeler89,Yu90}, which we now
describe.  Quasi-periodic extensions of general $t$-modules were
developed in \cite{BP02}.  Fix a Drinfeld $A$-module $\phi$ of rank
$r$ over $\Cinf$.  A $\phi$-\emph{biderivation} is an $\FF_q$-linear
map $\delta \colon A \rightarrow \CC_\infty\{\tau\}\tau$ satisfying
\[
\delta(ab) = a(\theta) \delta(b) + \delta(a)\phi(b),
\]
for all $a$, $b \in A$.  The space $D(\phi)$ of $\phi$-biderivations
splits as a direct sum
\[
D(\phi) = D_{sr}(\phi) \oplus \CC_\infty \delta_0 \oplus D_{si}(\phi),
\]
where
\begin{itemize}
\item $D_{sr}(\phi) = \{\delta \in D(\phi) \colon \deg_{\tau} \delta(a) < \deg_{\tau} \phi(a), \forall a \in A\}$,
\item $\delta_0(a) = \phi(a) - a\tau^0$,
\item $D_{si}(\phi) = \{\delta_T \in D(\phi)
\colon \delta_T(t) = T \phi(t) -\theta T,\ \text{some}\ T \in \Cinf\{\tau\}\tau\}$.
\end{itemize}
For each $\phi$-biderivation $\delta$, there is a unique entire
$\FF_q$-linear function $F_\delta(z)$, with no linear term, such
that
\[
F_\delta(\theta z) = \theta F_\delta(z) + \delta(t)\exp_\phi(z).
\]
The function $F_\delta$ is said to be the \emph{quasi-periodic} function related to $\delta$.
We note that  $F_{\delta_0}(z) = \exp_\phi(z) - z$ and that when
$\delta = \delta_T \in D_{si}(\phi)$, $F_{\delta_T}(z)$ is simply
$F_{\delta_T}(z) = T(\exp_\phi(z))$.

Corresponding to $\delta$ we define a $t$-module $\Phi_\delta$ by
\[
\Phi_\delta(t) =
\begin{pmatrix}
  \phi(t) & 0 \\
\delta(t) & \theta \tau^0
\end{pmatrix}
\]
with exponential function
\[
\Exp_{\Phi_\delta}
\begin{pmatrix}
  z \\ u
\end{pmatrix}
=
\begin{pmatrix}
  \exp_\phi(z) \\
u + F_\delta(z)
\end{pmatrix}.
\]
If $\omega$ is a period for $\exp_\phi(z)$, then
$(z,u) = (\omega, - F_\delta(\omega))$ is the corresponding period for
$\Exp_\delta$.  Therefore for $\phi = \phi_\Lambda$, the period lattice of $\Phi_\delta$ is
\[
\ker (\Exp_{\Phi_\delta}) =
\left\{
  \begin{pmatrix}
    \omega \\
-F_\delta(\omega)
  \end{pmatrix}\colon \omega \in \Lambda \right\}.
\]
In this way we can view $(\Phi_\delta,\GG_a^2)$ as an extension of $(\phi,\GG_a)$ by $\GG_a$:
\[
0 \rightarrow \GG_a \rightarrow (\Phi_\delta,\GG_a^2) \rightarrow (\phi,\GG_a) \rightarrow 0.
\]
This sequence splits in the category of $t$-modules precisely when
$\delta \in D_{in}(\phi) \as \Cinf \delta_0 \oplus D_{si}(\phi)$.

One can take these developments further and consider extensions of $\phi$ by several copies of $\GG_a$.  For example, we define the $t$-module
\[
\Phi(t) =
\begin{pmatrix}
  \phi(t) & 0 & 0 & \dots & 0 \\
\delta_1(t) & \theta \tau^0 & 0 & \dots & 0 \\
\delta_2(t) & 0 & \theta \tau^0 & \dots & 0 \\
\vdots & \vdots & & \ddots & \vdots \\
\delta_{s}(t) & 0 & \dots & 0 & \theta \tau^0  \\
\end{pmatrix},
\]
where $\delta_1, \dots, \delta_s \in D(\phi)$, which represents
an extension of $t$-modules
\[
 0 \to \GG_a^s \to (\Phi,\GG_a^{s+1}) \to (\phi,\GG_a) \to 0.
\]
Its corresponding exponential function is
\[
\Exp_\Phi \colon
\begin{pmatrix}
  z_0 \\ z_1 \\ \vdots \\ z_{s}
\end{pmatrix}
\longmapsto
\begin{pmatrix}
  \exp(z_0) \\ z_1+F_1(z_0) \\ \vdots \\ z_{s}+F_{s}(z_0)
\end{pmatrix}.
\]
One finds that, in the category of $t$-modules, $\Ext^1(\phi,\GG_a) \cong D(\phi)/D_{in}(\phi)$ as
$\Cinf$-vector spaces.  Gekeler~\cite{Gekeler89} found the dimension of
$D(\phi)/D_{in}(\phi)$ over $\Cinf$ to be $r-1$, and so any extension
$\Phi$ that contains no non-trivial subextensions must have $s
\leq r-1$.  For the connections among quasi-periodic functions and transcendence see \cite{Brown1998, Brown01, BP02, CP11, CP, Pellarin08, Yu90}.

\subsection{\texorpdfstring{$t$}{t}-modules from divided derivatives}
Until now examples of $t$-modules have been presented which in some
significant sense have analogues in the classical setting of abelian
varieties or commutative algebraic groups.  The present example is different in that it starts with a
$t$-module or a Drinfeld module and creates an extension in which the
coordinates of periods or logarithms of algebraic points are divided
derivatives of the coordinates of periods or logarithms of algebraic
points of the original $t$-module or Drinfeld module \cite{Brown01, BrDen00, Denis93}.

The \emph{divided derivatives} are a family $\{D_i\}$ of $\FF_q$-linear operators defined
first on $\FF_q[\theta]$ by the formula
\[
D_i \theta^n = \binom{n}{i} \theta^{n-i}.
\]
We extend $D_i$ to $k_{\infty}$ by continuity and then to $k_{\infty}^{sep}$.  Then
the divided derivatives satisfy the product formula
\[
D_n(ab) = \sum_{i+j = n} D_i(a) D_j(b).
\]
For legibility we often write $a^{[n]}$ for $D_n(a)$.

For simplicity of notation, we restrict ourselves to the case in
which we begin with a Drinfeld module $\phi$.  If $\phi(t) =
\sum_{h=0}^r a_h \tau^h$, then
the $t$-module $\Phi_{[n]}$, representing the divided derivatives of order up to $n$, is defined by
\[
\Phi_{[n]}(t) =
\begin{pmatrix}
  \phi(t) & 0 & \dots  &\dots & 0 \\
\phi^{[1]}(t) & \theta & & &\vdots  \\
\phi^{[2]}(t) & \Phi_{2,1}(t) & \ddots & & \vdots \\
 \vdots& \vdots& \ddots &\theta  & 0\\
\phi^{[n]}(t) & \Phi_{n,1}(t)   & \dots & \Phi_{n,n-1}(t) & \theta
\end{pmatrix},
\]
where (1) $\phi^{[i]}(t) = \sum_{h=0}^r a_h^{[i]} \tau^h$, (2) the matrix $\Phi_{[n]}(t)$
has zero superdiagonal terms, and (3) the $i+1,j+1$ terms below the main diagonal, with $1 \leq j < i \leq n$, are equal to
$\Phi_{i,j}(t) = \sum_{q^h \le i/j} a_{h}^{[i-q^h j]}\tau^h$. (One checks that the subdiagonal terms in the second through $n$-th columns are always $1$ unless $q=2$.)

The corresponding exponential function is
\[
\Exp_{[n]}
\begin{pmatrix}
  z_0 \\ z_1 \\ \vdots \\ z_n
\end{pmatrix}
=
\begin{pmatrix}
  \exp_\phi(z_0) \\
\exp_\phi^{[1]}(z_0) + z_1 \\
\vdots \\
\exp_\phi^{[n]}(z_0) + \sum_{q^h \le n}\sum_{r=1}^{\lfloor n/q^h\rfloor} a_h^{[n-q^hr]}z_r^{(h)}
\end{pmatrix}.
\]
If $\exp_\phi(u) \in k_\infty^{sep}$, then
\[
\Exp_{[n]}
\begin{pmatrix}
  u \\ u^{[1]} \\ \vdots \\ u^{[n]}
\end{pmatrix}
=
\begin{pmatrix}
\exp_\phi(u) \\
  \exp_\phi(u)^{[1]} \\
\vdots \\
\exp_\phi(u)^{[n]}
\end{pmatrix},
\]
using the fact that
\[
(x^q)^{[qi]} = (x^{[i]})^q.
\]
So if $\lambda$ is a period of $\exp_\phi$, then
\[
\begin{pmatrix}
  \lambda \\ \lambda^{[1]}\\ \vdots \\ \lambda^{[n]}
\end{pmatrix}
\]
is a period of $\Exp_{[n]}$.

\section{\texorpdfstring{$t$}{t}-Motives} \label{tMotives} %
Anderson gave definitions of two related, but different, kinds of
motives.  The first kind, called \emph{$t$-motives}, were introduced
in his seminal paper \cite{Anderson86}, and they have played a
dominant role in function field arithmetic.  The second, called
\emph{dual $t$-motives} have turned out to be more apt for recent
transcendence considerations \cite{ABP,ChangBanff,Papanikolas08}.  The reader
should consult \cite{Goss,Thakur} for additional information on $t$-motives.

\subsection{\texorpdfstring{$t$}{t}-motives}
Let $L$ be an extension of $\FF_q$, let $\iota \colon \FF_q[t] \to L$
be an $\FF_q$-algebra homomorphism, and set $\theta = \iota(t)$.
Let $L[t,\tau] \as L\{\tau\}[t]$ be the ring of polynomials in the
commuting variable $t$ over the non-commuting ring $L\{\tau\}$.  Thus
\[
  tc = ct, \quad t\tau = \tau  t, \quad \tau c = c^q \tau, \quad c \in L.
\]
A \emph{$t$-motive} $M$ is a left $L[t,\tau]$-module which is free and
finitely generated as an $L\{\tau\}$-module for which there is an
$\ell \in \NN$ with
\[
(t-\theta)^{\ell}(M/\tau M) = \{0\}.
\]
Strictly speaking, if $L$ is not perfect, we need to replace $\tau M$ by the $L[t]$-submodule it generates.
{\emph{Morphisms}} of $t$-motives are morphisms of left $L[t,\tau]$-modules.  The rank $d(M)$ of $M$ as an $L\{\tau\}$-module is called the \emph{dimension} of $M$.

Every $t$-module $(\Phi,\GG_{a}^{d})$ gives rise to the unique $t$-motive
\[
M(\Phi) \as \Hom^{q}_{L}(\GG_{a}^{d}, \GG_{a}) \cong L\{\tau\}^d,
\]
the module of $\FF_{q}$-linear morphisms  of algebraic groups.  The action of
 $L[t,\tau]$ is given by
\[
(ct^{i}, m) \longmapsto c \circ m \circ \Phi(t^{i}),
\]
$c \in L\{\tau\}$.  Projections onto the $d$ coordinates give an
$L\{\tau\}$-basis for $M$, $d = \rank_{L\{\tau\}} M$, and $\ell$
obviously need not be taken greater than $d$.

A $t$-motive has an $L\{\tau\}$-basis $m_{1},\dots, m_{d}$ which we can
use to express the $t$-action via a matrix $B(t) \in
\Mat_{d}(L\{\tau\})$.  This is compatible with the above considerations because, if
we represent an arbitrary element of $M$ as
\[
(k_{1},\dots,k_{d})
\begin{pmatrix}
  m_{1} \\ \vdots \\ m_{d}
\end{pmatrix} =
\mathbf{k}
\begin{pmatrix}
  m_{1} \\ \vdots \\ m_{d}
\end{pmatrix}, \quad k_i \in L\{\tau\},
\]
then according to the commutativity of $t$ with elements of
$L\{\tau\}$, for $c \in L\{\tau\}$,
\[
ct^{i}\cdot \mathbf{k}
\begin{pmatrix}
  m_{1} \\ \vdots \\ m_{d}
\end{pmatrix} =
c\mathbf{k}\cdot t^{i}
\begin{pmatrix}
  m_{1} \\ \vdots \\ m_{d}
\end{pmatrix} =
c\mathbf{k} B(t)^{i}
\begin{pmatrix}
  m_{1} \\ \vdots \\ m_{d}
\end{pmatrix}.
\]

A $t$-motive $M$ is said to be \emph{abelian} if it is free and finitely generated over $L[t]$.  In this case, the \emph{rank} of $M$ is its rank $r(M)$ as an $L[t]$-module.  A $t$-module is called abelian if its associated $t$-motive is abelian.

\begin{thm}[Equivalence of categories (Anderson
  {\cite[Thm.~1]{Anderson86}})] The above correspondence between abelian $t$-motives
  and abelian $t$-modules over $L$ gives an anti-equivalence of categories.
\end{thm}

An abelian $t$-motive $M$ is called \emph{pure} in the following situation.  We set
\[
 \laurent{M}{1/t} := M \otimes_{L[t]} \laurent{L}{1/t},
\]
which possesses a naturally induced left $L[t,\tau]$-module structure.  If $\laurent{M}{1/t}$ contains a finitely generated $\power{L}{1/t}$-submodule $H$, which generates $\laurent{M}{1/t}$ over $\laurent{L}{1/t}$ and which satisfies
\[
  t^u H = \tau^v H,
\]
for some $u$, $v \in \NN$, then $M$ is pure.  The \emph{weight} of $M$ is then
\[
  w(M) := \frac{d(M)}{r(M)} = \frac{u}{v}.
\]

When $L = \Cinf$, we can ask for criteria that govern when the $t$-module of an abelian $t$-motive is uniformizable.  Anderson shows that the notion of rigid analytic triviality for $t$-motives characterizes uniformizability.  For an integer $n$, we can define $n$-fold twisting of Laurent series in $\laurent{\Cinf}{t}$ by
\[
  f \mapsto f^{(n)} : \sum a_i t^i \mapsto \sum a_i^{(n)} t^i,
\]
and clearly each twisting map operates on both $\Cinf[t]$ and $\power{\Cinf}{t}$.  Moreover, the Tate algebra $\Cinf\langle t\rangle$ of all power series in $t$ over $\Cinf$ that converge on the closed unit disk in $\Cinf$ is also stable under twisting.
We can give $\laurent{\Cinf}{t}$ a ``trivial'' left $\Cinf[t,\tau]$-module structure by setting
\[
  \tau (f) := f^{(1)},
\]
and in this way $\Cinf[t]$, $\Cinf\langle t\rangle$, and $\power{\Cinf}{t}$ can also be given trivial left $\Cinf[t,\tau]$-module structures.  Now for an abelian $t$-motive $M$ we can make $M \otimes_{\Cinf[t]} \Cinf\langle t \rangle$ into a left $\Cinf[t,\tau]$-module by having $\tau$ act diagonally, using the above trivial $\tau$ action on $\Cinf \langle t \rangle$.  We say that $M$ is \emph{rigid analytically trivial} if
\[
  \bigl( M \otimes_{\Cinf[t]} \Cinf \langle t \rangle \bigr)^{\tau} \otimes_{\FF_q[t]}
   \Cinf \langle t \rangle \cong \bigl( \Cinf \langle t \rangle \bigr)^{r(M)},
\]
as left $\Cinf[t,\tau]$-modules.  Another way of putting this is as follows.  Let $\bm \in \Mat_{r \times 1}( \Cinf[t])$ be a basis for $M$ as a $\Cinf[t]$-module, and let $\Theta \in \Mat_r(\Cinf[t])$ represent multiplication by $\tau$ on $M$ with respect to $\bm$.  If $\Upsilon \in \GL_r(\Cinf\langle t \rangle)$ is chosen so that $\Upsilon \bm$ is a $\Cinf \langle t \rangle$-basis of $M \otimes_{\Cinf[t]} \Cinf \langle t \rangle$ fixed by the above action of $\tau$, i.e.\ $\tau(\Upsilon \bm) = \Upsilon \bm$, then we must have
\[
  \tau(\Upsilon \bm) = \Upsilon^{(1)} \Theta \bm = \Upsilon \bm,
\]
where the twisting $\Upsilon^{(1)}$ is applied entry-wise.  Thus $M$ is rigid analytically trivial precisely when there exists $\Upsilon \in \GL_r( \Cinf \langle t \rangle )$ that satisfies
\[
  \Upsilon = \Upsilon^{(1)} \Theta.
\]
We call such an $\Upsilon$ a \emph{rigid analytic trivialization} for $M$.
We then have the following fundamental result of Anderson.

\begin{thm}[Anderson {\cite[Thm.~4]{Anderson86}}]
The $t$-module associated to an abelian $t$-motive $M$ over $\Cinf$ is uniformizable if and only if $M$ is rigid analytically trivial.
\end{thm}

Now when $M_{1}$ and $M_{2}$ are pure $t$-motives, Anderson
constructs their \emph{tensor product} as the $t$-motive with
underlying module $M_{1}\otimes_{L[t]}M_{2}$ on which $\tau$ acts
diagonally.  Then $M_{1} \otimes_{L[t]} M_2$ is also a pure $t$-motive with weight
\[
  w(M_1 \otimes_{L[t]} M_2) = w(M_1) + w(M_2).
\]
In this way the category of pure $t$-motives over $L$ is a tensor
category.   It is not a Tannakian category (there is no trivial object nor are there
dual objects); however, it is possible to enlarge the category of
$t$-motives to a Tannakian category (see \cite{Papanikolas08, Taelman07}).  The construction of these categories can be reviewed in the current volume in~\cite[Rmks.~3.15, 4.14]{HartlJuschka}.

\subsection{The \texorpdfstring{$t$}{t}-motive of a Drinfeld module}
Let $\phi$ be a rank $r$ Drinfeld $A$-module defined over $L$ by
\[
  \phi(t) = \theta\tau^0 + a_1\tau + \dots + a_r\tau^r.
\]
Let $M(\phi) \as L\{\tau\}$, and as in the previous section we make $M(\phi)$ into the $t$-motive associated to $\phi$ by setting
\[
  ct^i \cdot m := cm\phi(t^i), \quad c \in L,\ m \in L\{\tau\}.
\]
We observe that $M(\phi)$ is an abelian $t$-motive.  Indeed we note that $1, \tau, \dots, \tau^{r-1}$ form an $L[t]$-basis for $M(\phi)$ (using the right division algorithm on $L\{\tau\}$).  In fact $M(\phi)$ is pure of dimension $1$, rank $r$, and weight $1/r$ \cite[Prop.~4.1.1]{Anderson86}.

\subsection{The \texorpdfstring{$t$}{t}-motive of \texorpdfstring{$C^{\otimes n}$}{C(n)}} \label{CarlitzTensorMot}
By the previous section the $t$-motive for the Carlitz module is $M(C) = L\{\tau\}$, where
\[
  t\cdot m = m(\theta + \tau), \quad m \in L\{\tau\}.
\]
Also $M(C)$ is rank $1$ over $L[t]$, generated by $1 \in L\{\tau\}$.  The function
\[
  \Omega(t) \as (-\theta)^{-q/(q-1)} \prod_{i=1}^\infty \left( 1 - \frac{t}{\theta^{(i)}} \right)
  \in \Cinf\langle t \rangle
\]
converges on all of $\Cinf$ and satisfies the difference equation
\[
  \Omega^{(-1)} = (t-\theta) \Omega.
\]
Thus $\Upsilon = (t-\theta)\Omega$ is a rigid analytic trivialization for $M(C)$.

The $t$-motive $M(C^{\otimes n})$ has rank $1$, dimension $n$, and weight $n$.  It is given by the $n$-fold tensor product
\[
  M(C^{\otimes n}) = L\{\tau\} \otimes_{L[t]} \cdots \otimes_{L[t]} L\{\tau\},
\]
on which $\tau$ acts diagonally.  It is rank $1$ over $L[t]$, generated by $1 \otimes \cdots \otimes 1$.  The calculation
\begin{align*}
  (t-\theta)(1 \otimes \cdots \otimes 1) &= 1 \otimes \cdots \otimes (t-\theta)\cdot 1 \otimes \cdots \otimes 1 \\
  &=1 \otimes \cdots \otimes \tau \otimes \cdots \otimes 1,
\end{align*}
where $t-\theta$ is multiplied at any arbitrary entry, implies that all elements of the form
\[
  1 \otimes \cdots \otimes \tau \otimes \cdots \otimes 1
\]
with one entry $\tau$ and the others $1$ are the same in $M(C^{\otimes n})$.  Repeating this construction by multiplying by additional factors of $t-\theta$ and applying an induction argument on the total degree of $\tau$'s appearing, we find that
\begin{align*}
m_1 &\as 1 \otimes \cdots \otimes 1 \\
m_2 &\as \tau \otimes 1 \otimes \cdots \otimes 1 \\
m_3 &\as \tau \otimes \tau \otimes 1 \otimes \cdots\otimes 1 \\
\vdots & \hspace*{30pt} \vdots \\
m_n &\as \tau \otimes \cdots \otimes \tau \otimes 1
\end{align*}
are a basis of $M(C^{\otimes n})$ as an $L\{\tau\}$-module.  The action of $t-\theta$ on this basis is
\begin{align*}
(t-\theta)m_i &= m_{i+1}, \qquad 1\leq i < n \\
(t-\theta)m_n &= \tau m_1.
\end{align*}
Thus
\[
  t \cdot \begin{pmatrix} m_1 \\ \vdots \\ m_n \end{pmatrix} =
  \begin{pmatrix} \theta & 1 & & \\ & \ddots & \ddots & \\ & & \ddots & 1 \\ \tau & & &\theta
  \end{pmatrix}
  \begin{pmatrix} m_1 \\ \vdots \\ m_n \end{pmatrix}
  = C^{\otimes n}(t) \begin{pmatrix} m_1 \\ \vdots \\ m_n \end{pmatrix},
\]
which agrees with the definition of $C^{\otimes n}$ in \S\ref{CarlitzTensorMod}.  Furthermore $(t-\theta)^n\Omega^n$ is a rigid analytic trivialization of $M(C^{\otimes n})$.

\subsection{Dual \texorpdfstring{$t$}{t}-motives}
In this section we define \emph{dual $t$-motives} and show how they are related to Drinfeld modules and $t$-modules.  In many ways the world of dual $t$-motives runs parallel (anti-parallel!)~with the world of $t$-motives already discussed.  There are some technical advantages, particularly for transcendence theory, for considering dual $t$-motives.  However, one disadvantage is that they are defined only over perfect fields. References for these types of $t$-motives are \cite[\S 4]{ABP} and \cite[\S 3.4]{Papanikolas08}.  The theory of dual $t$-motives was further developed by Anderson, who generously allowed this work to be reproduced in \cite[\S\S 5.1--5.2]{HartlJuschka}.

Let $L$ be an extension of $\FF_q$ that is perfect.  The ring $L[t,\sigma]$ is the polynomial ring in $t$ and $\sigma$ with coefficients in $L$ subject to the following relations,
\[
  tc = ct, \quad t\sigma = \sigma t, \quad \sigma c = c^{1/q}\sigma,
  \quad c \in L.
\]
In this way for any $f \in L[t]$,
\[
  \sigma f = f^{(-1)}\sigma.
\]
A \emph{dual $t$-motive} $H$ is a left $L[t,\sigma]$-module that is free and finitely generated over $L\{\sigma\}$ and for which there is an $\ell \in \NN$ with
\[
(t-\theta)^\ell (H/\sigma H) = \{0\}.
\]
Morphisms of dual $t$-motives are simply morphisms of left $L[t,\sigma]$-modules.  If in addition $H$ is free and finitely generated over $L[t]$, then $H$ is \emph{$A$-finite}.  Given a dual $t$-motive $H$, if the entries of $\bh \in
\Mat_{r \times 1}(M)$ comprise an $L[t]$-basis for $H$, then there is
a matrix $\Phi \in \Mat_r(L[t])$ so that
\[
  \sigma \bh = \Phi \bh.
\]
Since a power of $t-\theta$ annihilates $H/\sigma H$, we have
\[
  \det \Phi = c(t-\theta)^s,
\]
for some $c \in L^{\times}$ and $s \geq 0$.  In fact it can be shown that $s$ is the rank of $H$ as an $L\{\sigma\}$-module.

Just as in the case of $t$-motives we would like to describe explicitly the correspondence from $t$-modules to dual $t$-motives.  The map $* \colon L\{\tau\} \to L\{\sigma\}$ defined by
\[
  \left( \sum a_i \tau^i \right)^* \as \sum a_i^{(-i)} \sigma^i,
\]
is an anti-isomorphism of $L$-algebras.  For a matrix $B \in \Mat_{m \times n}(L\{\tau\})$ we define $B^* \in \Mat_{n\times m}(L\{\sigma\})$ by $(B^*)_{ij} := (B_{ji})^*$, and in this way
\[
  (BC)^* = C^* B^*
\]
whenever the matrix products are defined.  Now given a $t$-module $\Phi \colon A \to \Mat_d(L\{\tau\})$, $A = \FF_q[t]$, we let $H(\Phi) := \Mat_{1 \times d}(L\{\sigma\})$.  We make $H(\Phi)$ into a left $L[t]$-module by setting
\[
  a \cdot h := h\Phi(a)^*, \quad h \in H(\Phi),\ a \in A.
\]
One can check that $H(\Phi)$ is then a dual $t$-motive and that the evident functor $\Phi \mapsto H(\Phi)$ is covariant.  In fact Anderson has shown the following result.

\begin{thm}[Anderson]
  The correspondence between dual $t$-motives and $t$-modules over
  a perfect field $L$ gives an equivalence of
  categories.
\end{thm}

Recovering the $t$-module from its corresponding dual $t$-motive is quite easy.  In fact, if $H = H(\Phi)$, then the quotient module $H/(\sigma-1)H$ can be identified with $L^d$ as $\FF_q$-vector spaces.  If we examine the action of $t$ on this copy of $L^d$, we find that it acts by $\Phi(t)$, and so
\[
  \frac{H(\Phi)}{(\sigma -1)H(\Phi)} \cong (\Phi,L^d)
\]
as $A$-modules.  We will make this explicit in the upcoming examples.

\subsection{Rigid analytic triviality}
As in \S\ref{tMotives} there is also a corresponding notion of rigid analytic triviality for dual $t$-motives.  Working along closely similar lines we find that an $A$-finite dual $t$-motive $H$ over $\Cinf$, on which multiplication by $\sigma$ is represented by a matrix $\Phi \in \Mat_r(\Cinf[t])$, is \emph{rigid analytically trivial} if there exists $\Psi \in \GL_d(\Cinf\langle t \rangle)$ satisfying the system of difference equations
\begin{equation} \label{PsiPhi}
  \Psi^{(-1)} = \Phi\Psi.
\end{equation}
Here $\Psi^{(-1)}$ is the entry-wise twisting of $\Psi$.  And once again the $t$-module associated to $H$ is uniformizable precisely when $H$ is rigid analytically trivial.  Anderson's proof of this can be found in~\cite[Thm.~5.28]{HartlJuschka}.

As we will see in coming examples, when we evaluate $\Psi$ at $t=\theta$, we obtain a matrix whose entries are related to periods and quasi-periods of the $t$-module $E$ associated to $H$.  The chain of constructions take the following path:
\begin{align*}
  \Bigl\{ \textnormal{Uniformizable $A$-finite $t$-module $E$} \Bigr\}
  &\Longleftrightarrow
  \left\{ \parbox{2truein}{\begin{center} Rigid analytically trivial $A$-finite dual $t$-motive $H$ \end{center}} \right\} \\
  &\Longleftrightarrow
  \left\{ \parbox{2.5truein}{\begin{center} $\Phi \in \Mat_r(\Cinf[t])$
  representing $\sigma$ \\ with $\Psi \in \GL_r(\Cinf\langle t \rangle)$ satisfying\\
  $\Psi^{(-1)} = \Phi\Psi$ \end{center} }  \right\} \\
  &\Longrightarrow
  \left\{ \parbox{2truein}{\begin{center} $\Psi(\theta)^{-1}$ provides periods\\ and quasi-periods of $E$
  \end{center} }
  \right\}.
\end{align*}
We remark that we say $E$ is $A$-finite precisely when its dual $t$-motive is $A$-finite.  This notion is dual to the property of $E$ and its $t$-motive being abelian, but it is an open question whether $E$ being abelian and $E$ being $A$-finite are equivalent.  See the article by Hartl and Juschka \cite{HartlJuschka} in this volume for more details on these two properties.
The final implication, proved for periods by Anderson in unpublished work, allows us to arrive back at the $t$-module $E$ but now with precise information about its periods and quasi-periods.  For examples in several contexts, see \cite{ABP, BP02, CP11, CP, ChangYu07, Juschka, Papanikolas08, Pellarin08, Sinha97}.  These examples come from transcendence theory, as difference equations of the type in \eqref{PsiPhi} were studied extensively in \cite{ABP, Papanikolas08}.  See the article by C.-Y.\ Chang~\cite{ChangBanff} in this volume for additional information on these connections.

\subsection{Dual \texorpdfstring{$t$}{t}-motives for Drinfeld modules and \texorpdfstring{$t$}{t}-modules}
Suppose we have a Drinfeld $A$-module $\phi : A \to L\{\tau\}$, $A = \FF_q[t]$, given by
\[
  \phi(t) = \theta + a_1 \tau + \cdots + a_r \tau^r,
\]
and associated dual $t$-motive $H(\phi)$.  As in the previous section, we let $H(\phi) = L\{\sigma\}$ and then the action of $A$ on $H(\phi)$ is induced by
\[
  t \cdot b\sigma^j \as \sum_i a_i^{(-i-j)}b \sigma^{i+j}, \quad b \in L.
\]

By using this description we can define the dual $t$-motive
$H(C)$ that corresponds to the Carlitz module, and in this way for $h \in H(C) = L\{\sigma\}$,
\[
  t \cdot h = h (\theta+\sigma).
\]
Because $L$ is perfect, $L\{\sigma\}$ has a right division algorithm, and using it we can show that $1 \in L\{\sigma\}$ forms a basis for $H(C)$ as an $L[t]$-module.  From this point of view, given $h \in H(C)$, we can write
\[
  h  = f\cdot 1, \quad f \in L[t],
\]
and then
\[
  \sigma h = \sigma( f\cdot 1) = f^{(-1)} \cdot \sigma = f^{(-1)}(t-\theta)\cdot 1,
\]
since on $H(C)$ we have $t \cdot 1 = \theta + \sigma$.  Moreover, if we consider the image of $x \in L$ in $H(C)/(\sigma-1)H(C)$, we have the following calculation
\begin{align*}
  t \cdot x = \theta x + (t-\theta) \cdot x &= \theta x + \sigma(x^q)\\ &= \theta x +
  x^q + (\sigma -1)(x^q) \\ &= C(t)(x) + (\sigma-1)(x^q).
\end{align*}
Thus
\[
  \frac{H(C)}{(\sigma-1)H(C)} \cong (C,L)
\]
as $A$-modules.  We note that $\Omega(t)$ from \S\ref{CarlitzTensorMot} satisfies the difference equation
\[
  \Omega^{(-1)} = (t-\theta)\Omega,
\]
which makes $\Omega$ a rigid analytic trivialization of $H(C)$ in the dual $t$-motive framework (see \cite{ABP,Papanikolas08}).  We see further that by specializing at $t=\theta$, we obtain
\[
  \frac{1}{\Omega(\theta)} = -\tpi
\]
from \eqref{CarlitzPer}, which exemplifies the link between rigid analytic trivializations of $t$-motives and periods of their corresponding $t$-modules.

Now let $\phi$ be a rank $2$ Drinfeld $A$-module given by
\[
  \rho(t) = \theta + \kappa \tau + \tau^2, \quad \kappa \in L.
\]
Again using the functor from Drinfeld modules to dual $t$-motives, we arrive at the following construction.  The dual $t$-motive $H(\phi)$ is identified with $L\{\sigma\}$, on which $t$ acts by the rule
\[
  t\cdot h = h (\theta + \kappa^{(-1)}\sigma + \sigma^2), \quad h \in L\{\sigma\}.
\]
The right division algorithm on $L\{\sigma\}$ implies that $1$, $\sigma$ form an $L[t]$-basis of $H(\phi)$, and since
\[
  (t-\theta)\cdot 1 = \kappa^{(-1)}\sigma +\sigma^2,
\]
we have
\[
  \sigma \begin{pmatrix} 1 \\ \sigma \end{pmatrix} = \begin{pmatrix} \sigma \\ \sigma^2 \end{pmatrix}
  = \begin{pmatrix} 0 & 1 \\ t-\theta & -\kappa^{(-1)} \end{pmatrix}
  \cdot \begin{pmatrix} 1 \\ \sigma \end{pmatrix}.
\]
For $a$, $b \in L[t]$, we set
\[
  [a,b] := a\cdot 1 + b\cdot \sigma \in L\{\sigma\}.
\]
Now for $x \in L$,
\begin{align*}
  (\sigma -1)[x,0] &= -x + x^{(-1)}\sigma = \left[ -x, x^{(-1)} \right], \\
  (\sigma^2-1)[x,0] &= -x + x^{(-2)}\sigma^2 = \left[x^{(-2)}(t-\theta) -
  x,-\kappa^{(-1)}x^{(-2)}\right].
\end{align*}
Thus
\begin{align*}
  t[x,0] = [tx,0] &= \left[tx + \kappa x^q,-\kappa^{(-1)} x\right] + \left[-\kappa
  x^q,\kappa^{(-1)} x\right]\\
  &=\left[tx + \kappa x^q,-\kappa^{(-1)} x\right] + (\sigma - 1)[\kappa x^q,0]
  \\
  &= \left[\theta x + \kappa x^q + x^{q^2},0\right] + \left[(t-\theta)x -
  x^{q^2},-\kappa^{(-1)}x\right] + (\sigma - 1)\left[\kappa x^q,0\right] \\
  &=\left[\theta x + \kappa x^q + x^{q^2},0\right] + (\sigma - 1)\left[\kappa x^q,0\right]
  + (\sigma^2 -1)\left[x^{q^2},0\right].
\end{align*}
{From} this we see that the action of $t$ on $H(\phi)/(\sigma -1)H(\phi)$ is the same as the action of $\phi(t)$ on $L$, and we find that
\[
  \frac{H(\phi)}{(\sigma-1)H(\phi)} \cong (\phi,L)
\]
as $A$-modules.  Setting $\Phi = \left( \begin{smallmatrix} 0 & 1 \\ t-\theta & -\kappa^{(-1)} \end{smallmatrix} \right)$, the matrix that represents multiplication by $\sigma$ on $H(\phi)$, it is possible to find a matrix $\Psi \in \GL_2(\power{\Cinf}{t})$, whose entries are entire functions and which satisfies the difference equation
\[
  \Psi^{(-1)} = \Phi\Psi.
\]
This makes $\Psi$ into a rigid analytic trivialization for $H(\phi)$, and we find that
\[
  \Psi(\theta)^{-1} = \begin{pmatrix}
  \omega_1 & \eta_1 \\
  \omega_2 & \eta_2
\end{pmatrix},
\]
where $\omega_1$, $\omega_2$, $\eta_1$, $\eta_2$ are the fundamental
periods and quasi-periods of $\phi$ as in \S\ref{DrinfeldMods} and
\S\ref{QuasiPerMods}.  See \cite{CP11,Pellarin08} for complete
details.

\end{document}